\def\to{\rightarrow}
\def\La{{\Leftarrow}}
\def\Ra{{\Rightarrow}}

\def\be{\begin{equation}}

\def\sm{\setminus}
\def\ti{\tilde}
\def\pa{\partial}
\def\bs{\bigskip}\def\ss{\smallskip}\def\ms{\medskip}
\def\no{\noindent}

\def\PW{\text{\rm PW}}

\def\Res{\text{\rm Res}}
\def\const{\text{\rm const}}

\def\sign{\text{\rm sign}}

\def\ker{\text{\rm ker}}
\def\dim{\text{\rm dim}}

\def\dist{\text{\rm dist}}
\def\supp{\text{\rm supp}\,}

\def\AA{{\mathcal{A}}}

\def\DD{{\mathcal{D}}}
\def\EE{{\mathcal E}}

\def\HH{{\mathcal H}}

\def\MM{{\mathcal M}}
\def\NN{{\mathcal N}}

\def\SS{{\mathcal{S}}}

\def\C{{\mathbf{C}}}

\def\BM{{\mathcal BM}}

\def\D{{\mathbf{D}}}
\def\Z{{\mathbf{Z}}}

\def\R{{\mathbb R}}

\documentstyle[10pt,amssymb]{amsart}
\title{  
 Beurling-Malliavin theory for Toeplitz kernels }
\author{N.~Makarov}
\thanks{The first author is supported by 
N.S.F. Grant No. 0201893 }
\address{California Institute of Technology\\
Department of Mathematics\\
Pasadena, CA 91125, USA}
\email{makarov@@its.caltech.edu}
\author{A.~Poltoratski}
\address{Texas A\& M University
\\ Department of Mathematics\\
College Station, TX 77843, USA }
\email{alexeip@@math.tamu.edu}
\thanks{The second author is supported by 
N.S.F. Grant No. 0500852}

\theoremstyle{plain}
\newtheorem*{lem}{Lemma}
\newtheorem{le1}{Lemma}

\newtheorem*{thm}{Theorem}
\newtheorem*{thmA}{Theorem A}
\newtheorem*{thmC}{Theorem C}
\newtheorem*{thmB}{Theorem B}
\newtheorem*{cor}{Corollary}
\newtheorem*{prop}{Proposition}

\newenvironment{proof}
{{\noindent\it Proof:}}{\hfill$\Box$}

\numberwithin{equation}{section}

\begin{document}
\maketitle
\begin{abstract} We  consider the  family of Toeplitz operators $T_{J\bar S^{a}}$ acting 
in the Hardy space $H^2$ in the upper halfplane;  $J$ and $S$ are given meromorphic  inner functions, and $a$ is a real parameter.  In  the case where the argument of $S$ has  a power law type behavior on the real line,   
we compute the critical value 
$$ c(J,S)=\inf\left\{a:~ \ker~ T_{J\bar S^{a}}\ne0\right\}.$$
The formula for $ c(J,S)$ generalizes the Beurling-Malliavin theorem on  the  radius of completeness for a system of exponentials. 
\end{abstract}

\bs\section{Introduction and results}

\bs\subsection{Completeness of complex exponentials} For $\Lambda\subset\C$ denote
$$\EE_\Lambda=\left\{e^{i\lambda x}:~\lambda\in \Lambda\right\}.$$
By definition, the {\it radius of completeness} for the family $\EE_\Lambda$ is the number
$$R(\Lambda)=\sup\{a:~\EE_\Lambda\;{\rm is~complete~in}\; L^2(-a,a)\}.$$
(A family is complete if finite linear combinations of its elements are dense in the corresponding  space.)
In \cite{BM1}-\cite{BM2},  Beurling and Malliavin established a formula for  $R(\Lambda)$ in terms a certain density of $\Lambda$ at infinity.

\bs\no If $\Lambda\subset\R$, then the {\it effective} (or  Beurling-Malliavin)  density  $D_{\rm eff}(\Lambda)$ is  the supremum of the set of numbers $a\ge0$  such that there is a collection of disjoint intervals $\{l_j\}$ satisfying
the following two conditions:
$$\sum_j\frac{|l_j|^2}{1+d_j^2}=\infty,\qquad d_j:=\dist(0,l_j),$$
and 
$$\forall j,\qquad \#(\Lambda\cap l_j)\ge a |l_j|.$$

\ss\no 
Beurling-Malliavin's "{\it Second Theorem}" (BM2 for short) states that 
if $\Lambda\subset \R$, then $$R(\Lambda)=\pi D_{\rm eff}(\Lambda).$$

\ss\no 
 This formula extends to the general case $\Lambda\subset \C$ as follows. If $\Lambda$ satisfies the  {\it Blaschke} condition
$${\rm (B)}\qquad\qquad \sum_{\lambda\in\Lambda}\left|\Im~\lambda^{-1}\right|<\infty,$$then
$$R(\Lambda)=\pi D_{\rm eff}(\Lambda^*),$$where
$$\Lambda^*=\left\{\lambda^*:\;\lambda\in
\Lambda,\; \Re \lambda\ne0\right\},\qquad \lambda^*:=\left[\Re~\lambda^{-1}\right]^{-1};$$
otherwise $$\Lambda\not \in(B)\qquad\Ra\qquad R(\Lambda)=\infty.$$ 

\bs\no The Beurling-Malliavin theorem crowned a long search for a solution of the completeness problem, see \cite{PW}, \cite{L},  \cite{Sch}, \cite{Lv}. 
We refer to \cite{Re}  for  historical  information; let us only mention  that one of the earliest results of the theory was
the estimate  \be\label{udc}R(\Lambda)\le\pi D(\Lambda),\qquad (\Lambda\subset\R),\end{equation} where $D(\Lambda)$ is 
  the  usual {\it upper} density of $\Lambda$ at infinity.  

\bs\no The Beurling-Malliavin theory also comprises their "First Theorem" (BM1), a result of considerable independent interest and (so far) a necessary step in the proof of BM2. A detailed exposition of BM theory (including clarification and further improvements of the argument) is presented in the monographs \cite{Koo1},  \cite{Koo2}, 
\cite{HJ}, see also  \cite{Ka}, \cite{dB}. New applications and new approaches to various parts of the theory have  been suggested; see \cite{Br}, \cite{MH}, \cite {MNH} for some  recent developments; also see \cite{Kh} for a modern overview of the completeness problem for exponential systems. 

\bs\no In this paper we generalize BM theory to many other families of special functions. We state our results in the language of Toeplitz kernels referring to our paper \cite{MIF} for a detailed explanation of how  results of this type are related to the completeness problem for families of solutions of 
general Sturm-Liouville problems.

\bs\subsection{Toeplitz kernels}  The completeness radius problem can be restated in terms of Toeplitz operators as follows. Recall that the Toeplitz operator $T_U$ with   a symbol $U\in L^\infty(\R)$ 
is the map
$$T_U:H^2\to H^2,\qquad F\mapsto P_+(UF),$$
where $P_+$ is the orthogonal projection  in $L^2(\R)$ onto the Hardy space $H^2=H^2(\C_+)$ in the upper halfplane  $\C_+=\{\Im z>0\}$.
By  duality and the  definition of the
 classical Fourier transform, 
\begin{equation*}\label{cFF} f(t)\mapsto \hat f(z)=\int e^{izt}f(t)dt,\end{equation*} the exponential family   $\EE_\Lambda$ is complete in  $L^2(-a,a)$ if and only if there is a non-trivial function $F$ in the  Paley-Wiener space   
$$ \PW_a=\{\hat f:~f\in  L^2(-a,a)\}$$
such that $F=0$ on $\Lambda$.
 According to  Paley-Wiener's theorem, 
 the Fourier transform isometrically  identifies $L^2(0,\infty)$ with $H^2(\C_+)$, and therefore 
$$\PW_a=e^{-iaz}\left[H^2\ominus e^{2aiz}H^2\right].$$
The subspace $H^2\ominus e^{2aiz}H^2$ is the so called {\it model space} of the inner function $e^{2aiz}$. More generally, one defines  model spaces
$$K_\Theta=H^2\ominus \Theta H^2$$  for arbitrary inner functions $\Theta$. The elements of $K_\Theta$ are analytic functions in $\C_+$ but if $\Theta$ has a meromorphic extention to the whole complex plane, (we call such $\Theta$'s {\it meromorphic inner functions}), then the elements of  $K_\Theta$ are defined as   functions in $\C$. The completeness problem for exponentials is exactly the problem of describing the  sets of uniqueness for the model space of $e^{2aiz}$.

\ms\no 
Suppose now that $\Lambda$ is a subset of $\C_+$ satisfying the Blaschke condition, and let $B_\Lambda$ be the corresponding Blaschke
product. A simple argument shows that 
$\Lambda$ is a set of uniqueness for $K_\Theta$ if and only if  the Toeplitz operator with the symbol $U=B_\Lambda\bar \Theta $ has a trivial kernel. 
In particular, we obtain the formula 
 $$R(\Lambda)=\inf\left\{a:~ \ker~ T_{B_\Lambda e^{-2aiz}}\ne0\right\}.$$
There is a similar statement in the general case $\Lambda\subset\C$, see \cite{MIF}, Section 3.1. For example, if  $\Lambda\subset\R$, then 
$$R(\Lambda)=\inf\left\{a:~ \ker~ T_{J_\Lambda e^{-2aiz}}\ne0\right\},$$
where $J_\Lambda$ denotes some/any meromorphic inner function $J$ such that $\Lambda$ is precisely the level set $\{J=1\}$. 

\bs\no We should mention that the idea of the  Toeplitz operator approach in the study of exponential systems was introduced in the series of papers
 \cite{Pa}, \cite{Ni1}, \cite{HNP}. This approach has been particularly successful for the  interpolation and sampling theory in Paley-Wiener spaces,  see
\cite{LS}, \cite{OS}, \cite{Seip}.

\bs\no We will use the following notation for  kernels of  Toeplitz operators (or {\it Toeplitz kernels} in $H^2$):
$$N[U]=\ker~T_U.$$
(For example, $N[\bar\Theta]=K_\Theta$ if $\Theta$ is an inner function.)
We will also consider   Toeplitz kernels in the Smirnov-Nevanlinna  class $\NN^+=\NN^+(\C_+)$,
$$N^+[U]=\{F\in \NN^+\cap L^1_{{\rm loc}}(\R):~\bar U\bar F\in \NN^+\},$$
and in the  Hardy spaces $H^p=H^p(\C_+)$,
$$N^p[U]=N^+[U]\cap L^p(\R),\qquad (0<p\le\infty).$$
See \cite{Koo}, \cite{Ga}, \cite{Ni} for general    references concerning the   Hardy-Nevanlinna theory.

\bs\subsection{Generalization of Beurling-Malliavin theory}
\bs\no A natural way to generalize the completeness  radius problem (and the BM2 theorem) is  to  ask about the exact value of  the infimum
\be\label{g}\inf\left\{a:~ \ker~ T_{J\bar S^{a}}\ne0\right\}\end{equation}
for arbitrary meromorphic  inner functions $J$ and $S$. We will give an answer in the case where the argument of $S$ has  a power law type behavior,   
$$(\arg~ S)'(x)\asymp |x|^\kappa,\qquad  x\to\pm\infty,$$
with $\kappa\ge0$.  (We call the case $\kappa\ge 0$  {\it super-exponential} to underline  the relation to the classical case
$S(x)=e^{iax}$. In Section 1.5 below we will comment on the {\it sub-exponential} case $\kappa<0$.)

\ms\no As explained in  \cite{MIF}, the computation of the "radius" \eqref{g} has some immediate consequences for  the theory of Sturm-Liouville (SL) operators.  Roughly speaking, the case  of SL operators with  eigenvalues $$\lambda_n\asymp n^\nu$$ belongs to the theory with parameter
$$\kappa=\frac2\nu-1\ge0.$$
If $\kappa>0$, the SL operators are singular in contrast to the BM case $S(x)=e^{iax}$, which applies to regular  operators.
In addition to the completeness problem for systems of solutions of  SL equations, cf \cite{Hig}, the  generalized BM theory applies to certain  problems of  spectral theory as well as the theory of (Weyl-Titchmarsh) Fourier transforms associated with SL operators  and the corresponding (de Branges) spaces of entire functions.

\bs\no To state our results, we need to introduce the notion of {\it BM intervals}. ~Let $\gamma$ be a continuous  function $\R\to\R$ such that $\gamma(\mp\infty)=\pm\infty$. i.e.
$$\lim_{x\to-\infty}\gamma(x)=+\infty,\qquad \lim_{x\to+\infty}\gamma(x)=-\infty.$$ The family $\BM(\gamma)$ is defined as the collection  of the components of the open set
$$\left\{x:~\gamma(x)\ne\max_{[x,+\infty)}\gamma\right\}.$$ 
For an interval $l\in\BM(\gamma) $, we denote its length by $|l|$ or simply by  $l$, and we denote the distance to the origin by $d=d(l)$.

\bs\subsection{Super-exponential case}  If $\kappa\ge0$, then we say that $\gamma$ is $(\kappa)$-{\it almost decreasing} if \be\label{kap}\gamma(\mp\infty)=\pm\infty,\qquad \sum_{l\in\BM(\gamma)}d^{\kappa-2}l^2<\infty,\end{equation}
where the sum is taken over intervals satisfying $d(l)\ge1$. 

\ms\no The standard terminology in the classical $\kappa=0$ case is the following: the family
 $\BM(\gamma)$ is {\it short} if $\gamma$ is
almost decreasing; otherwise  we say that  $\BM(\gamma)$ is  {\it long}.

\ms
\begin{thmA} Let $\kappa\ge0$,  and let $U=e^{i\gamma}$ and  $S=e^{i\sigma}$ be smooth unimodular functions on $\R$ such that
\be\label{TA}\gamma'(x)\gtrsim -|x|^{\kappa},\quad \sigma'(x)\gtrsim|x|^{\kappa},\qquad (x\to\infty).\end{equation}

\ms\no (i) If $\gamma$ is not $(\kappa)$-almost  decreasing, then 
$N^+[US^\epsilon]=0$ for all $\epsilon>0$.

\ms\no (ii) If  $\gamma$ is  $(\kappa)$-almost  decreasing, then  $N^p[U\bar S^\epsilon]\ne0$ for all $\epsilon>0$   and all $p<\frac13$.
\end{thmA}

\ms\no Here and throughout the paper the notation $f(x)\gtrsim g(x)$  means that
$f(x)\ge cg(x)$ for some $c>0$ and all $x$ such that $|x|\gg1$.

\bs\no Given two unimodular functions $U$ and  $S$ as in Theorem A, we can consider the family of symbols$$U\bar S^a=e^{i\gamma_a},\qquad \gamma_a=\gamma-a\sigma,\qquad (a\in\R).$$
If $a_1>a$ and if $\gamma_a$ is decreasing near $\infty$, then $\gamma_{a_1}$ is also decreasing. It is not difficult  to see  that the same is true for almost decreasing functions, so  we can define the transition parameter 
$$c\equiv c(U,S;\kappa)=\inf\{a:~ \gamma_a\;{\rm is}~ \text{$(\kappa)$-almost descreasing}\}~\in (-\infty,+\infty].$$

\ms\begin{cor} Let $U=e^{i\gamma}$ and  $S=e^{i\sigma}$ be such that
$$\gamma'(x)\gtrsim -|x|^{\kappa},\quad \sigma'(x)\asymp|x|^{\kappa},\qquad (x\to \infty),$$
and let
$c=c(U,S;\kappa)$. Then for all $p<1/3$ we have
$$ N^p[U\bar S^a]=0 \quad (a<c),\qquad  N^p[U\bar S^a]\ne0 \quad (a>c).$$\end{cor} 

\bs\no Indeed, if $a<c$ then $\gamma_{a+\epsilon}$ is not almost decreasing for some $\epsilon>0$, and we have
$$ N^p[U\bar S^a]\subset N^+[U\bar S^{a+\epsilon}S^\epsilon]=0$$ by Theorem A, which can be applied because
 $\gamma_a'(x)\gtrsim -|x|^{\kappa}$ for all $a$'s. Similarly, if $a>c$, then $\gamma_{a-\epsilon}$ is  almost decreasing for some $\epsilon>0$, and we have
$$ N^p[U\bar S^a]=N^p[U\bar S^{a-\epsilon}\bar S^\epsilon]\ne0.$$

\bs\no In the special case where $U$ is an inner function, we can extend the  statement of he corollary to all values of $p$, in particular $p=2$.

\ss\begin{thmB} Let $J$ be a meromorphic  inner function, and suppose that a unimodular function $S$ satisfies $$(\arg~S)'(x)\asymp |x|^\kappa,\qquad x\to\infty.$$ Denote $c=c(J,S;\kappa)$. Then for all $p\le\infty$ we have
$$ N^p[J\bar S^a]=0 \quad (a<c),\qquad  N^p[J\bar S^a]\ne0 \quad (a>c).$$
\end{thmB}

\ms\subsection{Sub-exponential case} It is easy to see that the statement of Theorem A (and Theorem B) can not be extended to the case $\kappa<0$. For example,
the functions
$$\sigma(x)=2~\sign(x)~|x|^{1/4},\qquad \gamma(x)=2(1+\sqrt2)~ 1_{\R_-}(x)~|x|^{1/4}$$
satisfy the conditions \eqref{TA} with $\kappa=-3/4$, and of course $\gamma(+\infty)\ne-\infty$. Nevertheless, for $U=e^{i\gamma}$ and $S=e^{i\sigma}$ we have
$$N^\infty[US]\ne0,$$
 because
$$US=\bar f/f,\qquad f(z)=\exp\{-(1+i)z^{1/4}\}\in H^\infty(\C_+).$$
(Also note that the sum  
$\sum d^{\kappa-2}l^2$ in \eqref{kap} is always finite  if $\kappa<0$.)

\bs\no  The Beurling-Malliavin theory extends to  the sub-exponential case in a different fashion.
For $\kappa\in(-1,0]$ we consider the  weighted (non-linear) Smirnov-Nevanlinna classes 
$$\NN^+_\kappa=\left\{F\in \NN^+:~\log|F|\in L^1\left(\R,\frac1{1+|x|^{2+\kappa}}\right)\right\},$$
and  define the corresponding 
 Toeplitz kernels as follows:
$$N^+_\kappa[U]=N^+[U]\cap\NN^+_\kappa,\qquad N^p_\kappa[U]=N^p[U]\cap\NN^+_\kappa.$$

\ms\begin{thmC} Let $\kappa\in(-1,0]$,  and let $U=e^{i\gamma}$ and  $S=e^{i\sigma}$ be smooth unimodular functions  such that
$$\gamma'(x)\gtrsim -|x|^{\kappa},\quad \sigma'(x)\gtrsim|x|^{\kappa}\qquad (x\to \infty).$$

\ms\no (i) If the family $\BM(\gamma)$ is long, then 
$N^+_\kappa[US^\epsilon]=0$ for all $\epsilon>0$.

\ms\no (ii) If  the family $\BM(\gamma)$ is short, then  $N^p_\kappa[U\bar S^\epsilon]\ne0$ for all $\epsilon>0$   and all $p<\frac13$.
\end{thmC}

\ss\no  One can  also state a theorem similar to Theorem B.
Applications of these results to Volterra operators, see \cite{GK}, and higher order differential operators will be discussed in a separate paper.

\bs\subsection{Hilbert transform} The main tool in the proof of the theorems stated above is the one-dimensional  Hilbert transform.  Let $\Pi$ denote the Poisson measure on $\R$,
$$d\Pi(t)=\frac{dt}{1+t^2}.$$
If $h\in L^1_\Pi\equiv L^1(\R,\Pi)$ is a real-valued function, and if $\SS h$ denotes  its Schwarz integral,  
\be\label{schwarz}\SS h(z)=\frac1{\pi i}\int_\R\left[\frac1{t-z}-\frac t{1+t^2}\right]h(t)~dt,\qquad (z\in\C_+),\end{equation} 
then $\ti h$, the Hilbert transform  of $h$, is  defined a.e.~on $\R$ as the angular limit of $\Im [\SS h]$. 
Alternatively, $\ti h$ can be defined as a singular integral:
$$\ti h(x)=\frac1\pi~{\rm v.p.}\int\left[\frac1{x-t}+\frac t{1+t^2}\right]h(t)dt,\qquad (x\in\R).$$

\ss\no
(As a general rule we identify Nevanlinna class functions in the halfplane $\C_+$  with their angular boundary values on $\R$; e.g. we may write 
$\SS h=h+i\ti h$.)

\bs\no The relevance of the Hilbert transform in the theory of Toeplitz kernels can be explained by  the following simple observation, see \cite{MIF}, Section 2. 

\ms\no Suppose  $\gamma:\R\to\R$ is a smooth function. Then $N^+\left[e^{i\gamma}\right]\ne0$ if and only if 
\begin{equation}\label{basic} 
\gamma=-\alpha+\ti h\end{equation} for some smooth {\it increasing} function $\alpha$ and some  $h\in L^1_\Pi$. 
 There is a similar criterion for Toeplitz kernels in Hardy spaces: $N^p\left[e^{i\gamma}\right]\ne0$ 
if and only if  $\gamma$ admits a representation \eqref{basic} with $\alpha$ being the argument of some inner function and with $h\in L^1_\Pi$ such that $e^{-h}\in L^{p/2}(\R)$.

\bs\no 
 For further references, we recall some properties of the Hilbert transform.
We denote by $L^{(1,\infty)}_\Pi$ the usual weak $L^1$-space with respect to the Poisson measure.
 Kolmogorov's theorem states that 
$$\ti h\in L^{o(1,\infty)}_\Pi, $$ 
where $L^{o(1,\infty)}_\Pi$ stands for the "little o" subspace of $L^{(1,\infty)}_\Pi$,
i.e. 
\be\label{kol}\Pi\{|\ti h|>A\}=o\left(\frac1A\right),\qquad A\to\infty.\end{equation}
 For bounded functions we have the following (Smirnov-Kolmogorov) estimate : 
$$\|h\|_\infty<\frac\pi2\qquad\Ra\qquad e^{\ti h}\in L^1_\Pi.$$
Together with the criterion \eqref{basic}, this  implies
 \be\label{bmo2}\|\gamma\|_\infty<\frac\pi p\qquad\Ra\qquad N^p[\bar b^{2/p}e^{i\gamma}]\ne0,\end{equation}
 where $b$ is the Blaschke factor
\be\label{bla}b(z)=\frac{i-z}{i+z},\qquad z\in \C_+.\end{equation}

\bs\subsection{The structure of BM theory}  In the remaining sections of the paper we prove Theorems A and B.  We closely follow all the steps in our  presentation  of  the classical  Beurling-Malliavin theory in \cite{MIF}, combining them  with certain  estimates of the Hilbert transform, which we derive in Section~2. To make the proof self-contained, in several places we had to repeat the argument  outlined in \cite{MIF}. To  avoid further repetitions we decided to omit the proof of Theorem C because the  reasoning in the sub-exponential case is quite similar.
The proof of Theorems A and B is organized as follows.

\bs\no  (1)~ {\it Upper density estimate}:~  $\gamma(\pm\infty)\ne\mp\infty$ implies  $N^+[US^\epsilon]=0$;~ Section 3.1.\\ This statement is analogous to the estimate \eqref{udc}.

\bs\no  (2)~ {\it Effective density estimate}:~ $\sum d^{\kappa-2}l^2=\infty$ implies $N^+[US^\epsilon]=0$;~ Sections 3.2-3.4. Together with (1) this generalizes the estimate $R(\Lambda)\le \pi D_{\rm eff}(\Lambda)$ in BM2.

\bs\no  (3)~ {\it Little multiplier theorem}:~ if $\gamma$ is almost decreasing, then $N^+[U\bar S^\epsilon]=0$;~ Section~4. 

\bs\no(4)~ {\it BM multiplier theorem}:~ if the weighted Dirichlet norm of $\log W$ is finite, then $W$ belongs to some Hardy space up to a factor from $N^+[\bar S^\epsilon]$;~ Sections 5.1-5.2. 

\bs\no (5)~ {\it A version of BM1}:~ the logarithm of any outer function  in  $N^+$-kernel has a finite weighted Dirichlet norm. This is used to show that non-triviality of $N^+$-kernels implies 
 non-triviality of $N^p$-kernels for  symbols involving inner functions;~ Sections 5.3-5.4.
 
\bs \no(6)~ {\it $L^p$-multipliers}:~ approximation by inner functions and multiplying   the elements of $N^p$-kernels down to $H^\infty$;~ Section 6.

\bs
\section{One-sided Lipschitz condition for the  Hilbert transform}

\bs\no In this section we discuss various consequences of the weighted  one-sided Lipschitz condition
$$\ti h'(x)\lesssim |x|^{\kappa},\qquad x\to\infty,$$
for {\it smooth real-valued} functions  $h\in L^1_\Pi$.
(Here and elsewhere $\ti h'$ means $(\ti h)'$.) 

\bs\subsection{Application of Kolmogorov's theorem}
\begin{le1} If $\kappa\ge0$ and $h\in L_\Pi^1$, then 
$$\ti h'(x)\lesssim x^{\kappa}\quad \Rightarrow\quad \ti h(x)=o\left(x^{\kappa+1}\right)\qquad {\rm as}\quad x\to+\infty.$$\end{le1}

\ms
\begin{proof} By Kolmogorov's theorem, we have
\begin{equation}\label{Ko}\ti h\in L_\Pi^{o(1,\infty)}.\end{equation}  If $x_*\gg1$ and $\ti h(x_*)\ge cx_*^{\kappa+1}$ for some $c>0$,   then for all $x$ such that $1\ll x\le x_*$ we have 
$$\ti h(x)\ge \ti h(x_*)-a\int^{x_*}_xt^{\kappa}~dt\ge cx_*^{\kappa+1}-a(x_*^{\kappa+1}-x^{\kappa+1})
,$$ and it follows that 
$$\ti h(x)\gtrsim x_*^{\kappa+1}\ge x_*$$
 for all $x$ in some  interval $(x_{**},  x_*)$ of length $\asymp x_*$. Since $\Pi(x_{**},  x_*)\asymp1/x_*$, this contradicts \eqref{Ko}, see \eqref{kol}. 

\ms\no  If $\ti h(x_*)\le -cx_*^{\kappa+1}$ for some $c>0$, then by a similar argument we have
$$\ti h(x)\lesssim -x_*^{\kappa+1}\le x_*$$
 for all $x$ in some  interval $(x_*,  x_{**})$ of length $\asymp x_*$, which again contradicts   \eqref{Ko}.
\end{proof}

\bs\no We will also need the following version of this lemma.

\ms
\begin{le1} Let $h\in L_\Pi^1$, $\kappa\ge0$, and  $a\in\R$. If \begin{equation}\label{ko1}\ti h'(x)+ax^{-1}\ti h(x)\le x^{\kappa},\qquad x\gg1,\end{equation}
then
$$\ti h(x)=o\left(x^{\kappa+1}\right)\quad {\rm and}\quad \ti h'(x)\le x^{\kappa}+o\left(x^{\kappa}\right)\qquad {\rm as}\quad x\to+\infty.$$ 
\end{le1}

\begin{proof} Suppose we have $\ti h(x_*)\ge cx_*^{\kappa+1}$ for some $x_*\gg1$. Let $x_1$ be the smallest positive number such that  $\ti h(x_1)=cx_*^{\kappa+1}$, so we have $1\ll x_1\le x_*$ and $\ti h\le cx_*^{\kappa+1}$ on $(0,x_1)$. Together with \eqref{ko1}, this implies
$$\ti h'(x)\lesssim \frac{x_*^{\kappa+1}}{x_1},\qquad x\in \left(\frac{x_1}2,x_1\right).$$ Arguing as in the previous proof, we see that $\ti h\gtrsim x_*^{\kappa+1}\ge x_1$
 on  some  interval  of length $\asymp x_1$, which contradicts the weak $L^1$-estimate \eqref{kol}. 
The argument in the case $\ti h(x_*)\le -cx_*^{\kappa+1}$ is similar.  \end{proof}

\bs\subsection{Hilbert transform in weighted $L^1$-spaces}
\begin{le1}  If  $\kappa\in [-1,0)$ and  $h\in L^1\left({|x|^{-2-\kappa}}\right)$, then
$$\ti h'(x)\lesssim x^{\kappa}\quad \Rightarrow\quad \ti h(x)=o\left(x^{\kappa+1}\right),\qquad x\to+\infty.$$
\end{le1}

\begin{proof} If $\kappa\in (-1,0)$, then  
the weight $|x|^{-2-\kappa}$ satisfies the Muckenhoupt 
$(A_1)$ condition at infinity, and therefore we have
$$h\in L^1\left({|x|^{-2-\kappa}}\right) \quad\Ra\quad \ti h\in L^{o(1,\infty)}\left({|x|^{-2-\kappa}}\right),$$
see \cite{HMW}. (One can also give an elementary proof for this particular weight.)
We then argue as in the proof of Lemma 1. For example, 
 if $x_*\gg1$ and $\ti h(x_*)\ge cx_*^{\kappa+1}$,   then 
$\ti h\gtrsim x_*^{\kappa+1}$
 on some  interval $[x_{**},  x_{*}]$ of length $\asymp x_*$. The weighted length of this interval is $\asymp~x_*^{-1-\kappa}$, which contradicts the weak $L^1$-estimate.

\ms\no If $\kappa=-1$, then we consider the function
$$h_1(x)=xh(x)\in L^1_\Pi.$$
Since $\ti h_1(x)=x\ti  h(x)$, we have
$$\ti h_1'(x)=\ti h(x)+x\ti  h'(x)\le x^{-1}\ti h_1(x)+O(1),\qquad x\to+\infty.$$
By Lemma 2, we get $\ti h_1(x)=o(x)$ and therefore $\ti h(x)=o(1)$.
\end{proof}

\bs
\subsection{Persistence of 1-sided Lipschitz condition}

\begin{le1} Let  $f\in L^1_\Pi$, $0\not\in\supp~f$,  and let
$$0\le\alpha\le\beta,\qquad {\rm or}\qquad 0\le\beta<\alpha<2.$$ Denote
$$g(x)=|x|^{-\alpha}f(x).$$
Then
\begin{equation}\label{per1}\ti f'(x)\le (1+o(1))~|x|^{\beta}\quad\Ra\quad \ti g'(x)\le (1+o(1))~|x|^{\beta-\alpha},\end{equation}
and 
 \begin{equation}\label{per2}x^{-1}\ti f'(x)\le (1+o(1))~|x|^{\beta-1}\quad\Ra\quad  x^{-1}\ti g'(x)\le (1+o(1))~|x|^{\beta-\alpha-1}.\end{equation}
\end{le1}

\ss
\begin{proof} We will  prove \eqref{per1}  for $x\to+\infty$. The proof of the other cases is similar. Since the statement is trivial for $\alpha=0$, we will assume  $\alpha>0$.

 \ms\no It is clear that we can modify $f$ on any finite interval,  so we will assume that $f(x)=x^N$ near the origin for some $N\gg1$. If we specify $\ti f(0)=0$, then 
\be\label{per3}|x|^{-\alpha}\ti f(x)\in L^1_\Pi.\end{equation}Indeed, by Lemma 1 we have $\ti f=O(|x|^N)$, and therefore
$$|x|^{-\alpha}|\ti f|=|x|^{-\alpha}~|\ti f|^{\alpha/N}~|\ti f|^{1-(\alpha/N)}\lesssim |\ti f|^{1-(\alpha/N)}\in L^1_\Pi$$
 by Kolmogorov's estimate \eqref{kol}.

\ms\no 
Consider the  analytic function
$$u(z)+i\ti u(z):= z^{-\alpha}({f+i\ti f})(z), \qquad z\in \C^+,$$
where $z^{-\alpha}$ denotes the branch positive on $\R_+$. Note that
$$ u(x)=g(x), \quad \ti u(x)=|x|^{-\alpha}\ti f(x)\qquad {\rm for}\quad x\in\R_+,
$$
and 
$$u(x)=|x|^{-\alpha}~[f(x)\cos\alpha\pi+\ti f(x)\sin\alpha\pi]\qquad {\rm for}\quad x\in\R_-.$$
By \eqref{per3}, $$g-u\in L^1_\Pi,\qquad g-u=0\quad\rm  on \quad \R_+,$$
so if we define
$$\delta(x)=\ti g(x)- |x|^{-\alpha}\ti f(x),  $$
then
 $\delta=\ti g-\ti u$ on $\R_+$, and we have the following representation for the derivative:
$$\delta'(x)=\widetilde{(g-u)}'(x)=\int_{-\infty}^0\frac{c_1 f(t)+c_2 \ti f(t)}{(t-x)^2}~\frac{dt}{|t|^\alpha},\qquad (x>0).$$
By the dominated convergence theorem 
$$\delta'(x)=o(1),\qquad x\to +\infty,$$ in particular
 $$\delta'(x)=o\left(x^{\beta-\alpha}\right)\qquad\rm if \quad  \beta\ge\alpha.
$$ 
In the case $0\le\beta<\alpha<2$, we consider the integrals
involving $f$ and $\ti f$ separately. We have
\begin{align*}
\int_{-\infty}^{-1}\frac{|f(t)|}{(t-x)^2}~\frac{dt}{|t|^\alpha}
&\le\int_{-\infty}^{-x}\frac{1}{|t|^\alpha}~\frac{|f(t)|~dt}{|t|^2}+\frac1{x^\alpha}\int_{-x}^{-1} \frac{|t|^{2-\alpha}}{x^{2-\alpha}}~\frac{|f(t)|~dt}{|t|^2}
\\&\le \frac{1}{x^\alpha}\int_{-\infty}^{-x}\frac{|f(t)|~dt}{|t|^2}+\frac1{x^\alpha}\int_{-x}^{-1} \frac{|f(t)|~dt}{|t|^2}
=o\left(x^{-\alpha}\right)=o\left(x^{\beta-\alpha}\right).
\end{align*}
Since $\beta\ge0$, by Lemma 1 we have $$\ti f(t)=o\left(|t|^{1+\beta}\right),$$
and since
$$\int_{-\infty}^{-1}\frac{|t|^{1+\beta-\alpha}}{(t-x)^2}dt
\le\int_{-\infty}^{-x}|t|^{\beta-\alpha-1}dt+\frac1{x^2}\int_{-x}^{-1} 
|t|^{1+\beta-\alpha}dt\asymp x^{\beta-\alpha},$$
we have 
$$\int_{-\infty}^{-1}\frac{|\ti f(t)|}{(t-x)^2}~\frac{dt}{|t|^\alpha}=o\left(x^{\beta-\alpha}\right).$$
It follows that in all cases we have
$$\delta'(x)=o\left(x^{\beta-\alpha}\right), \qquad x\to +\infty,$$
and therefore
$$\ti g'(x)=x^{-\alpha}\ti f'(x)- \alpha x^{-\alpha-1}\ti f(x)+\delta'(x)
\le x^{\beta-\alpha}+o\left(x^{\beta-\alpha}\right).$$
\end{proof}

\bs
\subsection{A converse} We will also need a converse of \eqref{per1}. We  state it  only for the range of parameters that  will be   used later.

\begin{le1} Let $ w\in L^1_\Pi$, $0\not\in\supp w$,  and 
$$0<\alpha\le\beta, \qquad {\rm or}\qquad 1\le\alpha\le\min(2,\beta+1).$$ Denote
$$h(x)=|x|^{-\alpha}w(x).$$Then
$$\ti h'(x)\le (1+o(1))~|x|^{\beta-\alpha}\qquad\Ra\qquad \ti w'(x)\le (1+o(1))~|x|^{\beta}.$$\end{le1}

\bs\begin{proof} (a) The case $\beta\ge \alpha$.
Let $n$ be an even integer such that
$$\alpha_1:=n-\alpha\in[0,2).$$
Define  $$g(x):=x^{-n}{w(x)}=|x|^{-\alpha_1}h(x).$$
Since $\alpha_1<2$ and   $\beta_1:=\beta-\alpha\ge0$, we can apply Lemma 4 to $f=h$ and $g$ and  obtain the estimate
$$\ti g'(x)\le |x|^{\beta_1-\alpha_1}+\dots=|x|^{\beta-n}+\dots.$$
Since $\ti w(x)={x^n}\ti g(x)$,
we have
$$\ti w'(x)=nx^{n-1} \ti g(x)+x^{n} \ti g'(x)\le nx^{-1}\ti w(x)+|x|^{\beta}+\dots,$$
and by Lemma 2, 
$$\ti w'(x)\le(1+o(1))~|x|^{\beta}.$$

\bs\no (b) The case $\alpha\in [1,2]$ and $\beta-\alpha\in[-1,0]$. Note that this implies $\beta\ge0$.
Define the functions
$$g(x)= x^{-1}w(x), \qquad f(x)=xh(x),$$
so $$g(x)=|x|^{-\alpha_1} f(x),\qquad \alpha_1=2-\alpha\in [0,1].$$
Let us show that
\be\label{per5}x^{-1} \ti f'\le |x|^{\beta_1-1}+\dots,\qquad \beta_1:=\beta-\alpha+1.\end{equation}
 Since $$h\in L^1\left(\frac 1{|x|^{2-\alpha}}\right)\subset  L^1\left(\frac 1{|x|^{2+\kappa}}\right),\qquad 
\kappa:=\beta-\alpha,$$
 by Lemma 3 we have $$\ti h(x)=o(|x|^{\kappa+1}),$$ 
and since $\ti f(x)=x\ti h(x)$, we obtain \eqref{per5}:
$$x^{-1} \ti f'(x)=\ti h'(x)+x^{-1}\ti h(x)\le  |x|^{\kappa}+o(|x|^{\kappa}).$$
We can now apply Lemma 4 with parameters $\alpha_1$ and $\beta_1 $. (Note that  $f\in L^1_\Pi$ and the parameters  are admissible.)
By \eqref{per2} we get the estimate
$$x^{-1}\ti g'(x)\le |x|^{\beta_1-\alpha_1-1}+\dots=|x|^{\beta-2}+\dots,$$
and from  $\ti w(x)=x\ti g(x)$ we derive
$$\ti w'(x)=\ti g(x)+x\ti g'(x)\le x^{-1} \ti w(x)+|x|^\beta+\dots.$$
Applying Lemma 2 we conclude the proof.
\end{proof}

\bs
\section{Triviality of Toeplitz kernels}

\bs\no In this section we  prove the first part of Theorem~$A$, which gives a sufficient condition for the  triviality of a Toeplitz kernel. Let us fix $\kappa\ge0$ and  consider two unimodular functions
$U=e^{i\gamma}$ and $S=e^{i\sigma}$ on $\R~$  satisfying\be\label{reg}\gamma'(x)\gtrsim-|x|^\kappa,\qquad \sigma'(x)\gtrsim|x|^\kappa,\qquad (x\to\infty).\end{equation}

\bs\subsection{Upper density estimate}

\begin{prop}   If $N^+[US^\epsilon]\ne0$ for some $\epsilon>0$, then $\gamma(\mp\infty)=\pm\infty$. \end{prop}

\ms
\begin{proof} If $N^+[US^\epsilon]\ne0$, then by the basic criterion \eqref{basic} we have 
$$\gamma+\epsilon \sigma+\alpha=\ti h,\qquad \alpha'\ge0,\quad h\in L^1_\Pi.$$ 
Therefore,  $$\ti h'(x)\gtrsim \gamma'(x)\gtrsim-|x|^{\kappa},$$ and  $\ti h(x)=o\left(|x|^{\kappa+1}\right)$ by Lemma 1. It follows that
 $$\frac{\gamma(x)}x\lesssim \frac{\gamma(x)}x+\frac{\alpha(x)}x=-\epsilon \frac{\sigma(x)}x+o\left(|x|^{\kappa}\right)\lesssim - 
|x|^{\kappa},$$ which implies  $\gamma(\mp\infty)=\pm\infty$.
\end{proof}

\bs\subsection{Effective density estimate} 

\ms\no  Let $c>0$ be a fixed constant. 
For an interval $l\subset\R$ we denote by $l'$ and $l''$ the intervals of length $c|l|$ adjacent to $l$ from the left and from the right respectively, and we define
$$\Delta^*_l[\gamma]=\inf_{l''}\gamma-\sup_{l'}\gamma.$$

\ss\begin{lem}  Let $\epsilon >0$ and suppose 
\begin{equation}\label{bm}\gamma(\mp)=\pm\infty,\qquad \sum_{l\in\BM(\gamma)}d^{\kappa-2}l^2=\infty.\end{equation} Then there is a constant $c>0$ and there is a collection of disjoint  intervals $\{l_n\}$ in $[1,+\infty)$ or in $(-\infty,-1]$ such that
\begin{equation}\label{bm1}\sum d_n^{\kappa-2}l_n^2=\infty,\qquad 10 l_n\le d_n,\qquad \rm mult\{5l_n\}<\infty,\end{equation}
and 
\begin{equation}\label{bm2} \Delta^*_{l_n}[\arg(US^\epsilon)]\ge cd_n^{\kappa}l_n.\end{equation}
\end{lem}

  \ss\no Here $5l$ is the notation for the interval of length $5|l|$  concentric with $l$, and  mult$\{\cdot\}$ is the multiplicity of the covering.

\bs\begin{proof} Suppose the sum \eqref{bm} over BM intervals in $\R_+$ is infinite.  If there are infinitely many BM intervals $l=(\tilde a_n, b_n)$ in $\R_+$   satisfying $10|l|> d$, then we set $$l_n=(a_n, b_n), \qquad a_n:=\frac{10}{11}~b_n;$$otherwise we simply enumerate BM intervals   such that $10|l|\le d$. In any case, we get a collection of intervals $l_n=(a_n,b_n)$ satisfying  the first two conditions in \eqref{bm1} and also the inequality
$$\gamma(b_n)\ge \gamma(a_n).$$
By \eqref{reg}, the latter implies  that the intervals also satisfy \eqref{bm2} for some $c>0$. Finally, we take a subfamily $\{l_{n_k}\}$ such that $\{5l_{n_k}\}$ is a subcover of $\bigcup 5l_n$ of finite multiplicity and observe that  we still have the divergence of the series
$\sum d^{\kappa-2}l^2$. Indeed, if~ $\bigcup l_j\subset 5l$,~ then  $$d_j\asymp d,\qquad \sum l_j^2\lesssim l^2,$$ and  so
$$\sum_j d_j^{\kappa-2}l_j^2\lesssim d^{\kappa-2}l^2.$$\end{proof}

\bs\no The following proposition completes the proof of the first part of Theorem A. 
\begin{prop} Suppose  $\gamma'(x)\gtrsim-|x|^\kappa$ and suppose there is a collection $\{l\}$ of disjoint intervals  in $[1,+\infty)$ such that   $$\forall l,\qquad\Delta^*_l[\gamma]\ge cld^{\kappa},$$ and \begin{equation*}\sum d^{\kappa-2}l^2=\infty,\qquad 10 l\le d,\qquad \rm mult\{5l\}<\infty,\end{equation*}
then $N^+\left[e^{i\gamma}\right]=0.$\end{prop}

\bs\subsection{Proof of the proposition}
The statement corresponds to the so-called Beurling's lemma in the classical BM theory. There are several  versions of the  proof of Beurling's lemma, e.g.    Koosis \cite{Koo1} applies the Beurling-Tsuji estimate of harmonic measure, Nazarov \cite{N}  uses the Bellman function, and Kargaev's proof \cite{Kar}  is  based on PDE techniques. We suggest yet another approach.

\bs\no According to the   criterion \eqref{basic}, we have to exclude the possibility
$$\gamma+\alpha=\tilde h,\qquad \alpha~\uparrow, \quad h\in L^1_\Pi.$$
 Denote by $ h_l$ the restriction of $h$ to the interval $5l$. We say that $l$ is {\it of type} I if \begin{equation}\label{type1} d^{\kappa-2}l^2\le C\|h_l\|_\Pi,\end{equation}
where $C$ is a sufficiently large constant; otherwise we call $l$ an interval {\it of type} II.
Clearly, we have 
\begin{equation*}\sum_{l\in\rm I} d^{\kappa-2}l^2<\infty,\end{equation*}
and to get a contradiction we need  to show\begin{equation}\label{type}\sum_{l\in\rm II} d^{\kappa-2}l^2<\infty\end{equation}

\ms\no Consider the 2D Hilbert transform 
$$H(z)=\int_\R\frac{h^-(t)~dt}{(t-z)^2},\qquad (z\in\C_+).$$
where $h^-=\max\{0,-h\}$.

\begin{lem} If $l$ is of type II, then $$|H(z)|\gtrsim d^{\kappa},\qquad \forall z\in Q_l:=\{x+iy: x\in l, \;l<y<2l\}.$$\end{lem}

\ss\no We prove this lemma in the next subsection, and we now explain how the lemma implies \eqref{type}.
Denote $$\psi=\sum_{l\in \rm II}d^{\kappa} l\cdot 1_l.$$
We have
\begin{equation}\label{83}\sum_{l\in\rm II} d^{\kappa-2} l^2\asymp \int_1^\infty\frac{\psi(t)~dt}{t^2}=\frac83\int_1^\infty\frac{dA}{A^3}\int_{A/2}^A\psi(t)~dt.\end{equation}For every $A>1$ let
$$H_A(z)=\int_{-CA}^{CA}\frac {h^-(t)~dt}{(t-z)^2},$$
where $C>0$ is a large constant, and let II$(A)$ denote the set of all intervals $l\in \rm II$ intersecting $(A/2,A)$. If $l\in {\rm II}(A)$ and $z\in Q_l$, then
$$|H(z)-H_A(z)|\le \int_{|t|>CA}\frac {h^-(t)~dt}{|t-z|^2}\asymp  \int_{|t|>CA}\frac {h^-(t)~dt}{1+t^2}\ll1,$$
so by the lemma we have
$$|H_A|\gtrsim A^{\kappa}\qquad{\rm on} \quad \bigcup_{l\in{\rm II}(A)} Q_l.$$
Applying the   weak-$L^1$ estimate for the 2D Hilbert transform, see \cite{CZ},  we get
$$\sum_{l\in{\rm II}(A)}l^2\le{\rm Area}(|H_A|\gtrsim A^{\kappa})\lesssim A^{-\kappa}\int_{-CA}^{CA}h^-(t)~dt,$$
and therefore
$$\int_{A/2}^A\psi(t)~dt\lesssim A^{\kappa}\sum_{l\in{\rm II}(A)}l^2
\lesssim\int_{-CA}^{CA}h^-(t)~dt,$$
Combining this with \eqref{83}, we conclude
$$\sum_{l\in\rm II} d^{\kappa-2} l^2\lesssim\int_1^\infty\frac{dA}{A^3}\int_{-CA}^{CA}h^-(t)~dt\lesssim \|h^-\|_\Pi,$$
which proves \eqref{type}

\bs\subsection{Proof of the lemma} Since $\alpha$ in the representation $\tilde h=\gamma+\alpha$ is increasing, we have
\begin{equation}\label{D1}\Delta^*_l[\ti h]\gtrsim ld^{\kappa}.\end{equation}
On the other hand, for intervals  of type II we have
\begin{equation}\label{D2}\Delta^*_l[\ti h_l]\ll ld^{\kappa}.\end{equation}
Indeed, if $~\Delta^*_l[\ti h_l]\gtrsim ld^{\kappa}~$, then $~|\ti h_l|\gtrsim 
 ld^{\kappa}~$ on either $l'$ or $l''$. Applying the weak type inequality with $A\asymp ld^{\kappa}$, we get
$$d^{-2}l~\lesssim ~\Pi\{|\ti h_l|>A\}~\lesssim~ A^{-1}\|h_l\|_\Pi,$$
which contradicts the definition of type II.

\ms\no 
 Denote $f\equiv f_l=1_{\R\sm 5l}\cdot h$, so 
$\ti h=\ti h_l+\ti f_l.$ From \eqref{D1}-\eqref{D2} we conclude that
 there are points $a\in l'$ and $b\in l''$ such that
\begin{equation}\label{l1}\ti f(b)-\ti f(a)\ge \frac c2 ld^{\kappa}.\end{equation}
Represent $f=f^+-f^-$ with $f^+=\max \{f,0\}$, and note that  the functions\\ $\tilde f_\pm:=(f^\pm)\ti{~}$
are decreasing on $[a,b]$:
$$\tilde f_\pm'(x)=-\frac1\pi\int_{\R\sm(5l)}\frac{f^\pm(t)~dt}{(t-x)^2}<0,\qquad (x\in 5l).$$
From \eqref{l1} it then follows that 
\begin{equation*}\ti f_-(a)-\ti f_-(b)\gtrsim ld^{\kappa},\end{equation*}
so there is a point $x_*\in (a,b)$
such that
$$\frac1\pi\int\frac{f^-(t)dt}{(t-x_*)^2}=-\tilde f_-'(x_*)=\frac{\ti f_-(a)-\ti f_-(b)}{b-a}\gtrsim d^{\kappa}.$$
Observe that if $z\in Q_l$ and $t\in\R\sm5l$, then
$$\Re\left[\frac 1{(t-z)^2}\right]\asymp \frac 1{(t-x_*)^2},$$
and we have
$$\left|\int\frac{f^-(t)dt}{(t-z)^2}\right|\ge\Re~\int\frac{f^-(t)dt}{(t-z)^2} \asymp\int\frac{f^-(t)dt}{(t-x_*)^2}
\gtrsim d^{\kappa}.$$
It follows that
$$|H(z)|\ge\left|\int\frac{f_l^-(t)dt}{(t-z)^2}\right|-\left|\int\frac{h_{l}^-(t)dt}{(t-z)^2}\right|  \gtrsim d^{\kappa},$$
because
$$\left|\int\frac{h_l^-(t)dt}{(t-z)^2}\right|~\lesssim\frac {d^2}{l^2}\|h_l\|_\Pi\ll d^{\kappa}$$
provided that the constant $C$ in \eqref{type1} is large enough.

\bs\section{Non-triviality of Toeplitz kernels in Smirnov-Nevanlinna  class}

\bs\subsection{A version of the "little multiplier" theorem} In this section we  prove the following statement. Let $\kappa\ge0$  and suppose that $U=e^{i\gamma}$, $S=e^{i\sigma}$ satisfy conditions \eqref{TA} of Theorem A.

\ss
\begin{prop} If $\gamma$ is almost $(\kappa)$ decreasing, then 
$N^+[U\bar S^\epsilon]\ne0$ for all $\epsilon>0$.\end{prop}

\ss\begin{proof} By assumption we have \begin{equation}\label{sum}\sum_{l\in\BM(\gamma)} d^{\kappa-2}l^2<\infty.\end{equation}
Recall that the   BM intervals of $\gamma$ are the components of the open set
$\{\gamma^*\ne \gamma\}$, where $\gamma^*(x)=\max\gamma[x,+\infty)$. 
Denote $$f=\gamma^*-\gamma,$$ so $f=0$ outside the union of BM intervals, and $f'(x)\lesssim |x|^{\kappa}$ on BM intervals.
By \eqref{sum} we have $l\lesssim d$, and therefore  
\begin{equation}\label{e1}0\le f\lesssim ld^{\kappa}\qquad{\rm on}\quad l,\end{equation}
Together with \eqref{sum} this implies   $$f\in L^1_\Pi.$$ 
The  estimate \eqref{e1} also shows that we can assume $l\ge d^{-\kappa}$ for  all BM intervals; otherwise we can eliminate short intervals   by adding a bounded function to $\gamma$ (this will not affect the $N^+$-kernel). In particular, we will assume that   BM intervals don't cluster to a finite point.

\ss\no The non-triviality of $N^+[U\bar S^\epsilon]$ is a consequence of the following statement which will be verified  in the next two subsections.

\ss
\begin{lem} {\it For any $\epsilon>0$, there is a
function $\beta$ such that} $$f+\beta\in\ti L_\Pi^1,\qquad \beta'(x)\le \epsilon |x|^{\kappa}\quad{\rm for}\quad|x|\gg1. $$\end{lem}

\ss\no Indeed, if for instance $\sigma'(x)\ge |x|^{\kappa}$ near $\pm\infty$, then we can write
$$\gamma-\epsilon\sigma=-(f+\beta)+(\beta-\epsilon\sigma)+\gamma^*.$$
The first term in the RHS is in $\ti L_\Pi^1$, and the last two terms are decreasing near infinities, so we can apply the basic criterion \eqref{basic}.\end{proof}

\bs\subsection{Proof of the lemma} 
\ms\no 
We will construct disjoint intervals $l_n$ such that they cover all BM intervals and satisfy the following two conditions:
\begin{equation}\label{e2}\sum_{n}d^{\kappa-2}_n l_n^2<\infty,\end{equation}
and
\begin{equation}\label{e3}\forall n\quad\exists \epsilon_n\in[0,\epsilon],\qquad
\int_{ l_n}\frac{f(x)-\epsilon_n |x|^{\kappa}T_{n}(x)}{1+x^2}~dx=0,\end{equation} 
where $T_n$ is the "tent" function of the interval $l_n$,
$$T_n(x)=\dist (x, \R\sm l_n).$$

\ms \no
Let us show that the existence
of such intervals $l_n$ implies the statement of the  lemma.
 Define
$$\beta(x)=-\sum_n\epsilon_n |x|^{\kappa}T_{n}(x),$$ and $$ g(x)=f(x)-\sum\epsilon_n |x|^{\kappa}T_{n}(x).$$
Clearly, we have
$$|\beta'(x)|\lesssim \epsilon|x|^{\kappa},$$
and all we need is to check   $ g\in \ti L^1_\Pi$. 

\ms\no Let us  show that $g$ belongs to the {\it real} Hardy space $\HH^1_\Pi(\R)$. We can represent $g$ as follows:
$$g=\sum g_n=
\sum\lambda_n\frac{g_n}{\lambda_n}:=\sum\lambda_nA_n,$$
where  
$$g_n=g\cdot1_{l_n},\qquad \lambda_n=\Pi(l_n)~\|g_n\|_\infty.$$
The functions $A_n=\lambda_n^{-1}g_n$ are "atoms":
$$\int A_n~d\Pi=\frac1{\lambda_n} \int_{ l_n}g~d\Pi=0\qquad {\rm by} \quad\eqref{e3},$$
and 
$$\|A_n\|_\infty=\frac{\|g_n\|_\infty}{\lambda_n}=\frac1{   \Pi(l_n)}.$$
Since $$\|g_n\|_\infty\lesssim d_n^{\kappa} l_n,$$
(use $l_n\lesssim d_n$ and \eqref{e1} for the BM intervals covered by $l_n$),
we  have
$$\sum\lambda_n \lesssim\sum\frac
{l_n}{d_n^2}~d_n^{\kappa}l_n<\infty\qquad {\rm by} \quad\eqref{e2}.$$
It follows that $\sum\lambda_nA_n\in \HH^1_\Pi(\R)$, see \cite{CW}.

\bs\subsection{Construction of  intervals $l_n$} Let us assume that
all BM intervals $l$ lie in $[1,+\infty)$. In the general case we will need to apply the procedure described below to BM intervals in $(-\infty,-1]$ and in $[1,+\infty)$ separately. 

\ms\no We construct our intervals $l_1, l_2,\dots$ by induction. The left endpoint $a_1$ of $l_1$ will be the left endpoint of the leftmost BM interval. 
 Suppose the left endpoint $a_n$ of $ l_n$ has been constructed so that $a_n$ is also the left endpoint
of some  BM interval $l=(a_{(l)}, b_{(l)})$, i.e. $a_n=a_{(l)}$.
 Consider the function
$$F(b)=\int_{a_n}^b\left[f-\epsilon |x|^{\kappa}T_{(a_n,b)}\right]~d\Pi,$$
where  $T_{(a_n,b)}(\cdot)=\dist(~\cdot~, \{a_n,b\})$ is the tent function.  
We 
define $b_n$, the right endpoint of $l_n$, as the nearest point in the complement of BM intervals at which $F$ is non-positive, 
$$ b_n=\min \{b\ge b_{(l)}:~f(b)=0,\;F(b)\le0\}.$$
 Since $f\in L^1_\Pi$, we have $F(+\infty)=-\infty$ and so $b_n<\infty$. 
Finally, we define   $a_{n+1}$ as the leftmost endpoint of BM intervals not covered by $l_1\cup \dots\cup l_n$. (Recall that we assumed that there are no finite cluster points.)

\bs\no It is clear from the construction that the intervals $l_n$ cover all BM intervals. We also get \eqref{e3} by defining  $\epsilon_n$ from the equation
$$\int_{a_n}^{b_n}\left[f-\epsilon_n |x|^{\kappa}T_{(a_n,b_n)}\right]~d\Pi=0;$$
clearly we have $0<\epsilon_n\le\epsilon$.
In remains to verify \eqref{e2}. We have  three types of intervals $l_n$:

\ms\no (a) $F( b_n)<0$ but there is a BM interval $l\subset l_n$ such that 
$|l|\asymp |l_n|$,

\ms\no (b) $F( b_n)=0$,

\ms\no (c) other intervals.

\ms\no Property \eqref{e2}  is obvious for the collection of   intervals of type (a): we have $l\ll d$ (except for finitely many $l$'s) and therefore $d\asymp d_n$ and $d^{\kappa-2}l^2\asymp d_n^{\kappa-2}l_n^2$.

\ms\no  To prove \eqref{e2} for the collection of  intervals of type (b), we 
note that $\epsilon_n=\epsilon$ if $l_n\in$(b), and since$$\int_a^bx^{\kappa} T_{(a,b)}(x)~\frac{dx}{x^2}\ge \frac{a^{\kappa}}{b^2}\int_a^b T_{(a,b)}\gtrsim \frac{a^{\kappa}}{b^2}(b-a)^2,$$
we have
 
$$\sum_{(b)}\frac{d_n^{\kappa}l_n^2}{(d_n+ l_n)^2}\lesssim \sum_{(b)}\int_{a_n}^{b_n}|x|^{\kappa}T_{(a_n,b_n)}~d\Pi=
\frac1\epsilon\int_{\bigcup_{(b)} l_n}f~d\Pi<\infty.$$
Since $d_n\to\infty$, it follows that there are only finitely many intervals $l_n\in(b)$ satisfying $d_n\le l_n$, so the last estimate implies $$\sum_{(b)}d_n^{\kappa-2}l_n^2<\infty.$$

\ss\no The argument for intervals of type (c) is the same if we can show that  if $l_n\in$(c), then $\epsilon_n>\epsilon/2$, i.e.
\begin{equation}\label{cp}\int_{a_n}^{b_n}\left[f-\frac\epsilon 2 |x|^{\kappa}T_{(a_n,b_n)}\right]~d\Pi>0.\end{equation} Since $l_n$ is not of type (b), we have
$F( b_n)<0$ and  by construction, $b_n$ is the right endpoint of some BM interval $l=(c,b_n)$. Note that $|l|\ll |l_n|$ because $l_n$ is not of type (a).
Since $f>0$ on $l$, we have
\begin{align*}
\int_{a_n}^{b_n}&\left(f- \frac\epsilon 2  x^{\kappa}
T_{(a_n,b_n)}\right)~d\Pi>\int_{a_n}^{c}\left(f-\epsilon x^{\kappa}
T_{(a_n,c)}\right)~d\Pi+ \\
&+\left[\epsilon\int_{a_n}^{c}x^{\kappa}
T_{(a_n,c)}~d\Pi
-\frac\epsilon 2\int_{a_n}^{b_n}x^{\kappa}
T_{(a_n,b_n)}~d\Pi
\right].\end{align*}
The first term in the RHS is equal to $F(c)$ and therefore positive by construction. Since $|l|\ll |l_n|$,  the second term in the RHS is also positive, and we get \eqref{cp}

\bs\section{Multiplier theorem}

\bs\subsection{BM multipliers} Let $S$ be a unimodular function and let  $0<p\le\infty$. If  $w\in L^1_\Pi$ is a real function, then we write
$$w\in \MM_p(S)$$ if the outer function $$W=e^{w+i\ti w}$$ satisfies the following condition:$$\quad\forall \epsilon>0,\qquad \exists G\in N^+[\bar S^\epsilon],\qquad WG\in H^p.$$
In other words, $w\in \MM_p(S)$ if the corresponding outer function belongs  to $H^p$ up to an arbitrarily small (compared to $S$) factor.

\ms\no 
 We can restate  this property  in terms of Toeplitz kernels.

\begin{lem} $w\in \MM_p(S)$ iff 
$$\forall \epsilon>0,\qquad N^p\left[\bar S^\epsilon \frac{\bar W}W\right]\ne0.$$\end{lem}
\begin{proof} $\Ra$\; Let $G\in N^+[\bar S^\epsilon]$ be such that
$F:=GW\in H^p$. Then
$$F\in N^p\left[\bar S^\epsilon \frac{\bar W}W\right],$$
and the Toeplitz kernel is non-trivial. Indeed,
$$\bar S^\epsilon \frac{\bar W}WF=(\bar S^\epsilon G)\bar W\in \NN^-\cap L^p=\bar H^p.$$

\ss\no  $\La$\; If $F$ is in the Toeplitz kernel, i.e. $F\in H^p$ and $F\bar S^\epsilon \bar W/W\in \bar H^p$, then we define $G=F/W\in \NN^+$. Since 
$$ \bar S^\epsilon G=\bar S^\epsilon\frac{\bar W}W F\frac1{\bar W}\in \NN^-,$$
we have $G\in N^+[\bar S^\epsilon]$ and $WG\in H^p$.
\end{proof}

\ms\begin{cor} Suppose  $(\arg S)'\gtrsim |x|^\kappa$. If a real function $w_0\in L^1_\Pi$ satisfies the following condition: $$\forall \epsilon>0, \quad\exists w\in  L^1_\Pi,\quad w\ge w_0, \quad \ti w'>-\epsilon|x|^{\kappa}+o\left(|x|^{\kappa}\right),$$
then $w_0\in \MM_p(S)$ for all $p<1$.
\end{cor}

\begin{proof} Without loss of generality, $(\arg S)'\ge 2|x|^\kappa$ for $|x|\gg1$. We have
$$-\arg~\left[\bar S^{2\epsilon} \frac{\bar W_0}{W_0}\right]=2\left(\epsilon\arg S+\ti w\right)+2(\ti w_0-\ti w):=\alpha+\ti g,$$
where $$\alpha'\ge \epsilon|x|^\kappa,\qquad(|x|\gg1),$$
and $$g\in L^1_\Pi,\qquad g\le0.$$ 
For sufficiently large $N$, the function $\alpha+\arg b^N$, where $b$ is the Blaschke factor \eqref{bla}, is monotone increasing on $\R$, and therefore there is an inner function $\Phi$, not  a finite Blaschke product, such that
$$\alpha+\arg b^N=\arg\Phi+\delta,\qquad \|\delta\|_\infty\le\pi.$$
Clearly,  $N^\infty\left[e^{-i\ti g}\right]\ne0$, i.e.
$$N^\infty\left[e^{i\delta}~\bar b^N~ \Phi~ \bar S^{2\epsilon}~ \frac{\bar W_0}{W_0}\right]~\ne~0.$$
By \eqref{bmo2}, we also have 
$$N^p\left[e^{-i\delta}~\bar b^N \right]\ne0,$$
provided that $p<1$ and $N$ is sufficiently large, and of course
$$N^\infty\left[b^{2N}~ \bar\Phi\right]\ne0. $$
It follows that $$N^p\left[\bar S^{2\epsilon}~ \frac{\bar W_0}{W_0}\right]~\ne~0.$$
\end{proof}

\ms\no 
The main result of this section is the following version of the Beurling-Malliavin multiplier theorem.

\ms
\begin{thm}  Suppose  $(\arg S)'\gtrsim |x|^\kappa$, and let $w_0\in L^1_\Pi$ be a real function. Then  
\begin{equation*} |x|^{-\frac{2+\kappa}2}w_0(x)\in\DD(\R,\infty)\qquad\Ra\qquad w_0\in \MM_p(S),\quad (\forall p<1).\end{equation*}
\end{thm}

\ms\no
Here the notation $f\in \DD(\R,\infty)$ means that there is a neighborhood of infinity where  $f$ coincides with some function  from the {\it Dirichlet} space $\DD(\R)$. Recall that the Hilbert space $\DD(\R)$ consists of functions $h\in L^1_\Pi$ such that the harmonic extension $u=u(z)$ of $h$ to $\C_+$ has a finite gradient norm, 
$$\|h\|^2_\DD\equiv \|u\|^2_\nabla=\int_{\C_+} |\nabla u|^2~dA<\infty,$$
($dA$ is the area measure). If  $h\in\DD(\R)$ is a smooth function, then we also have
$$\|h\|^2_\DD=\int_\R \bar h~\ti h'~dx.$$
In the next two subsections we use some ideas from the proof of Theorem 64 in \cite{dB}.

\bs\subsection{Proof of the multiplier theorem}   
It is clear that  we can assume that the function  $w_0+i\ti w_0$ is  analytic   and  has a zero of sufficiently large multiplicity  at the origin; in particular 
$$h_0(x):=|x|^{-\frac{2+\kappa}2}w_0(x)\in\DD(\R).$$
Let us fix $\epsilon>0$. According to the last corollary we need to construct $w$ such that
$$(i)\; w\in L^1_\Pi,\qquad (ii)\; w\ge w_0, \qquad (iii)\; \ti w'>-\epsilon |x|^{\kappa}+o\left(|x|^{\kappa}\right).$$
We define
$$w(x)=|x|^{\frac{2+\kappa}2}h(x),$$
where $h$ is a solution of the following extremal problem:
$$\min\{I(h):~h\ge h_0\},\qquad I(h):=\|h\|^2_{\DD}+\epsilon
\int|x|^{\frac{2+\kappa}2}|h(x)|~d\Pi(x)
.$$
The existence of a solution follows from the usual argument: the set $$\AA=\{h:~\|h\|_\DD\le I(h_0),\quad h\ge h_0 \;{\rm a.e.}\}\subset \DD(\R)$$
 is bounded, closed, and convex in $\DD(\R)$, therefore it is weakly compact. Let $I_0$ denote the minimum of $I(h)$ over $\AA$. Then there is a sequence of functions $h_n\in\AA$ such that $I(h_n)\to I_0$ and $h_n$ weakly converge to some function $g\in\AA$. It is then routine to see that
$$I_0\le I(g)\le \liminf I(h_n)=I_0,$$
so $g$ is a solution of the extremal problem.

\bs\no By construction, $w$ satisfies  (i) and (ii). To prove (iii) we first note that
\begin{equation}\label{bmm} \ti h'(x)\ge -\epsilon|x|^{\frac{\kappa-2}2}.\end{equation}
Indeed, 
by the  extremality of $h$ we have
$$ \|\phi\|^2_\DD+2\int \phi\ti h'+\epsilon
\int(|h+\phi|-|h|)(x)~\frac{|x|^{\frac{2+\kappa}2}}{1+x^2}~dx=I(h+\phi)-I(h)\ge0$$
for all smooth test functions $\phi=\phi(x)\ge0$. (The integral $\int \phi\ti h'$  has to be interpreted in the sense of the theory of distributions.) Since 
$$\frac{\phi(x)}{x^2}\ge\frac{|h(x)+\phi(x)|-|h(x)|}{1+x^2},$$
we conclude
$$\|\phi\|^2_\DD+2\int \phi(x)~\left[\ti h'(x)+\epsilon
|x|^{\frac{\kappa-2}2}\right]~dx\ge0$$
Replacing $\phi(x)$ with $\delta\phi(x)$ and letting $\delta\to0$, we get
$$\int \phi(x)~\left[\ti h'(x)+\epsilon
|x|^{\frac{\kappa-2}2}\right]~dx\ge0$$
for all $\phi\ge0$, which proves \eqref{bmm}

\bs\no To derive (iii) from \eqref{bmm} we apply Lemma 5 in Section 2 with
$$\alpha=1+\frac\kappa2,\qquad \beta=\kappa, \qquad \beta-\alpha=\frac\kappa2-1.$$
The parameters $\alpha$ and $\beta$ are admissible because for  $\kappa\ge2$ we have $\alpha\le \beta$, and if $0\le\kappa\le2$ then 
$1\le\alpha\le2$ and  $\alpha\le \beta+1$.

\bs\subsection{Multipliers and one-sided Lipschitz condition}

\ms
\begin{prop}  If  $w\in L^1_\Pi$, $w\ge0$, and $\ti w'\lesssim |x|^{\kappa}$, then 
$$|x|^{-\frac{2+\kappa}2}w(x)\in\DD(\R,\infty).$$
\end{prop}

\ms
\begin{proof} We will assume that the function  $w_0+i\ti w_0$ is  analytic   and  has a zero of sufficiently large multiplicity  at the origin. Let $u=u(z)$ be the harmonic extention of $|x|^{-\frac{2+\kappa}2}w(x)$ to the upper half plane $\C_+$, and let $v=\ti u$. We need to show that the gradient norm of $u+iv$ in  $\C_+$ is finite,
$$\|u+iv\|_{\nabla}^2=\lim_{r\to\infty}\int_{\pa D(r)}udv<\infty,$$
where $D(r)$ is the semidisc $\{|z|<r\}\cap\C_+$. 

\ms\no 
We first prove that the integrals over $\pa D(r)\cap\R$ 
are uniformly bounded from above. Applying Lemma 4 in Section 2 with  (admissible) parameters$$\alpha=1+\frac{\kappa}2,\qquad \beta=\kappa,\qquad \beta-\alpha=\frac{\kappa}2-1$$ to the functions $f=w$ and $g=u$, we see that
\begin{equation*}\label{1bm}v'(x)\lesssim |x|^{\frac{\kappa-2}2},\qquad x\in\R.\end{equation*}
Since $u>0$ we  have
$$\int_{\pa D(r)\cap\R}udv=\int_{-r}^r  v'u\lesssim\int_{-r}^r |x|^{\frac{\kappa-2}2}~|x|^{-\frac{2+\kappa}2}~w(x)~dx\lesssim \|w\|_\Pi
<\infty.$$

\bs\no To finish the proof of the proposition it remains to show that the integrals
$$\int_{\pa D(r)\sm\R}u~dv=rI'(r),\qquad
I(r):=\frac12\int_0^\pi u^2\left(re^{i\theta}\right)d\theta,$$
don't  tend to $+\infty$ as $r\to\infty$.
In fact, it is enough to show 
$$I(r)\not\to\infty,$$ because if~ $rI'(r)\to +\infty$, then $I'(r)\ge1/r$ for all  $r\gg1$, and we have $I(r)\to\infty$. Since $\kappa\ge0$, we can apply the following lemma. 

\bs\no{\bf Lemma.} {\it If $u\in L^1(1+|x|^{-1})$, then}  $I(r)\not\to\infty$.

\bs\no
{\it Proof:}~ We will prove an equivalent statement for functions in the unit disc $\D$. Let $f=u+i\ti u$ be an analytic function in $\D$ such that $$\frac {u(\zeta)}{1-|\zeta|}\in L^1(\pa\D).$$
Define
$$h(z)=\frac {1+z}{1-z}~u(z),\qquad z\in \D,$$
and denote by $h^*(\zeta)$, $\zeta\in\pa\D$, the angular maximal function. By Hardy-Littlewood maximal theorem,
\be\label{w*}h^*\in L^1_{\rm weak}(\pa\D).\end{equation} Let us show that as $\epsilon\to0$, 
$$\frac1\epsilon\int_{C_\epsilon}|f(z)|^2|dz|\not\to\infty,\qquad C_\epsilon:=\{|1-z|=\epsilon\}\cap\D.$$
We have
$$\frac1\epsilon\int_{C_\epsilon}|f|^2=\epsilon\int_{C_\epsilon}|h|^2\lesssim\left[\epsilon h^*(\zeta)\right]^2+\left[\epsilon h^*(\bar\zeta)\right]^2,$$
where $\zeta\in\pa\D$, $|1-\zeta|=\epsilon$.
The RHS can not tend to infinity because otherwise for all small $\epsilon$, we would have 
$$h^*(\zeta)+h^*(\bar\zeta)\gg\frac 1\epsilon$$
on an interval of length $\epsilon$, which would contradict \eqref{w*}.
\end{proof}

\bs\subsection{A version of BM1} 

\ms\begin{prop} Suppose    $(\arg S)'\gtrsim|x|^\kappa$ and let $\Theta$ be a meromorphic  inner function satisfying $|\Theta'|\lesssim|x|^{\kappa}$. Then$$W\in N^+[\bar\Theta]\qquad \Ra\qquad \log|W|\in \MM_p(S),\quad  (\forall p<1).$$ \end{prop}

\ms
\begin{proof} We have
$W\bar\Theta=\bar H$ for some $H\in\NN^+$. Define
$$W_1=WH+\Theta,$$
and let $W_1^e$ be the outer part of $W_1$.
From the identity
$$\bar\Theta^2W_1=\bar\Theta W\bar\Theta H+\bar\Theta=\bar H\bar W+\bar\Theta=\bar W_1,$$
we deduce 
$$|W_1|=|W\bar W\Theta+\Theta|=1+|W|^2\ge1,$$
and  $$|W_1^e|\ge1,\qquad |W|\le|W_1^e|,\qquad (\arg W_1^e)'\lesssim|x|^{\kappa}.$$
By Proposition 5.3 and the multiplier theorem, we have $\log|W_1|\in\MM_p(S)$ and therefore $$\log|W|\in
\MM_p(S).$$\end{proof}

\ms\begin{cor} Let $S$ and $\Theta$ be as above. Then for any
meromorphic inner function $J$ and any $p<1$, we have
$$N^+[\bar \Theta J]\ne0\qquad\Ra\qquad \forall \epsilon,\quad N^p[\bar S^\epsilon\bar \Theta J]\ne0.$$
\end{cor}

\begin{proof} Take  an outer function $W\in N^+[\bar\Theta J]$. 
Then $W\in N^+[\bar\Theta]$,
and by the last proposition, 
$$\exists G\in N^+[\bar S^\epsilon].\qquad WG\in H^p.$$
It then follows that
$$WG\in N^+[\bar S^\epsilon\bar\Theta J]\cap H^p=N^p[\bar S^\epsilon \bar\Theta J].$$
\end{proof}

\bs
\section{Non-triviality of Toeplitz kernels in Hardy spaces}

\bs\no In this final section we finish the proof of Theorems A and B.

\bs\subsection{Approximation by inner functions} 
It is well known that 
given any two intertwining discrete sets $A=\{a_n\}$ and $B=\{b_n\}$ of real numbers, \linebreak $...a_n<b_n<a_{n+1}...$,
there exists a meromorphic inner function $\Theta$ such that \be\label{AB}\{\Theta=1\}=A,\qquad \{\Theta=-1\}=B.\end{equation}
Indeed, the sequences $A$, $B$ determine the set
$$E=\{\Im\Theta>0\}\cap\R=\cup (a_n,b_n),$$
and we can define $\Theta$ in $\C_+$  by the (Krein's shift) formula
\be\label{ks}\frac 1{\pi i}\log\frac{\Theta+1}{\Theta-1}=\SS u+ic, \qquad u:=1_E-\frac12,\quad c\in\R,\end{equation}
where $\SS u$ is the Schwarz integral \eqref{schwarz}, so $\Re [\SS u]$ is the Poisson extension of $u$ to the halfplane. 
(Note that $u$ is the boundary function of the expression in the LHS of \eqref{ks}, provided that $\Theta$ is an inner function with level sets $A$ and  $B$, and  in fact Krein's shift  formula parametrizes all such inner functions.)

\ms\no 
An immediate consequence of this construction is the following statement:

\ss\no  {\it for any 
increasing, continuous  function $\sigma:\R\to\R$, 
there exists a meromorphic inner function $\Theta=e^{i\theta}$ such that }
$$\|\theta-\sigma\|_\infty\le\pi.$$

\ss\no  We will need the following  version of this statement. 

\begin{lem} If $\sigma'(x)\asymp |x|^\kappa$,
then there is a meromorphic inner function $\Theta=e^{i\theta}$ such that 
$$\theta-\sigma\in L^\infty(\R),\qquad \theta'(x)\asymp|x|^\kappa.$$\end{lem}

\begin{proof} We can assume that $\sigma$ is strictly increasing on $\R$. Define the  intertwining sequences $A=\{a_n\}$ and $B=\{b_n\}$   by the equations
$$\sigma(a_n)=2\pi n, \qquad b_n=\frac{a_n+a_{n+1}}2,\qquad (n\in\Z),$$
so we have 
$$a_n\asymp (\sign~ n)~|n|^{\frac1{1+\kappa}},$$
and 
$$ \delta_n:=b_n-a_n\asymp |a_n|^{-\kappa}.$$
 Let $\Theta$ be an inner function satisfying \eqref{AB}, $$\|\theta-\sigma\|_\infty\le2\pi,$$
and let $\mu_1$,  $\mu_{-1}$ be the corresponding  (Aleksandrov-Clark's) measures defined by the  Herglotz representation
$$\frac{1+\Theta}{1-\Theta}=\SS\mu_1+\const ,\qquad 
\frac{1-\Theta}{1+\Theta}=\SS\mu_{-1}+\const.$$
The measures $\mu_1$,  $\mu_{-1}$ have the following form:
$$\mu_1=\sum\alpha_n\delta_{a_n},\qquad \mu_{-1}=\sum\beta_n\delta_{b_n}$$ for some positive numbers $\alpha_n$, $\beta_n$. (It is easy to see that $\mu_{\pm 1}\{\infty\}=0$ though we don't actually need this fact.) We claim that
\be\label{ab} \alpha_n\asymp \delta_n, \qquad \beta_n\asymp \delta_n.\end{equation}
 The estimate 
$\theta'(x)\asymp|x|^\kappa$ easily follows from \eqref{ab}. 
 Since 
$$|\Theta'|\asymp \left|1-\Theta^2\right|~\left|(\SS\mu_1)'\right|,\qquad 
|\Theta'|\asymp \left|1+\Theta^2\right|~\left|(\SS\mu_{-1})'\right|,$$
we have
$$\theta'(x)\asymp\min \left\{ \sum\frac{\alpha_n}{(x-a_n)^2},\;  \sum\frac{\beta_n}{(x-b_n)^2}\right\},\qquad (x\in\R).$$
It follows that if $x\in(a_m, a_{m+1})$, then by \eqref{ab}
$$\theta'(x)\asymp
\int_{|t-x|\gtrsim\delta_m}\frac{dt}{(x-t)^2}\asymp \delta_m^{-1}\asymp |a_m|^\kappa\asymp |x|^\kappa.$$

\bs\no {\it Proof of \eqref{ab}.} We will explain  the estimate for $\alpha_n$'s;  the proof for $\beta_n$'s is similar. According to \eqref{ks}, we have
$$\frac{1-\Theta}{1+\Theta}=\const~ e^{Ku},$$
where $u=1_E-1/2$,$$E=\bigcup_{k=-\infty}^\infty~(a_k, b_k),$$ and $Ku$ is the {\it improper} integral 
$$Ku(z)=\int\frac {u(t)~dt}{t-z}, \qquad(z\in\C_+).$$
 By construction,   $$\alpha_n=\const~\Res_{a_n}e^{Ku}.$$ Denote 
$$g_n(z)=\exp\left\{\int_{b_{n-1}}^{b_{n}}\frac {u(t)~dt}{t-z}\right\}= 
\frac{\sqrt{(b_n-z)(b_{n-1}-z)}}{a_n-z},$$
and $$A_n=\exp\left\{\int_{\R\sm(b_{n-1},b_{n})}\frac {u(t)~dt}{t-a_n}\right\},$$
so
$$\Res_{a_n}e^{Ku}=A_n~\Res_{a_n}g_n, \qquad \left|\Res_{a_n}g_n\right| \asymp \delta_n.$$
It remains to show that 
$A_n=e^{O(1)}$. This can be done as follows.

\ms\no For $j>n$ we  have
\begin{align*}\int_{a_j}^{a_{j+1}}\frac {u(t)~dt}{t-a_n}&=\log\frac{b_j-a_n}{a_j-a_n}-\log\frac{a_{j+1}-a_n}{b_j-a_n}\\&=
\log\left(1+\frac{\delta_j}{a_j-a_n}\right)-\log\left(1+\frac{\delta_j}{b_j-a_n}\right)\\&=
\frac{\delta_j}{a_j-a_n}-\frac{\delta_j}{b_j-a_n}+O\left(\frac{\delta_j^2}{(a_j-a_n)^2}\right)=O\left(\frac{\delta_j^2}{(a_j-a_n)^2}\right),\end{align*}
where we used the relation $\log(1+x)=x+O(x^2)$ for $0<x\lesssim 1$. Since 
\begin{align*}\sum_{j=n+1}^\infty\frac{\delta_j^2}{(a_j-a_n)^2}&\asymp \sum_{j=n+1}^\infty\frac{\delta_j}{a_j^\kappa(a_j-a_n)^2}\\&\asymp \int_{a_n+\delta_n}^\infty \frac{dt}{t^\kappa(t-a_n)^2}=\int_{a_n+\delta_n}^{2a_n}
+\int_{2a_n}^\infty \\&\lesssim \frac1{a_n^\kappa}\int_{a_n+\delta_n}^{\infty}\frac{dt}{(t-a_n)^2}+\int_{a_n}^\infty\frac{dt}{t^{2+\kappa}}\\&\asymp 
\frac1{a_n^\kappa}~\frac1{\delta_n}+\frac1{a_n^{1+\kappa}}=O(1),
\end{align*}
we get$$\int_{b_{n}}^\infty\frac {u(t)~dt}{t-a_n}=O(1).$$
A similar estimate holds for the integral over $(-\infty, b_{n-1})$, and we have $A_n=e^{O(1)}$.
\end{proof}

\bs\subsection{Proof of Theorem A} The first part of the theorem was established in  Section 3. The second part states that if $\gamma$ is almost decreasing and $\epsilon >0$, then 
\begin{equation}\label{H1}N^p[U\bar S^{2\epsilon}]\ne0,\qquad (p<1/3).\end{equation} By Lemma 6.1 there exists  an inner function  $\Theta$ satisfying
$$(\arg \Theta)'\asymp |x|^\kappa.$$
We will assume that $U^2\Theta$ has an increasing argument (otherwise we can replace $\Theta$ with $\Theta^n$ for a large  integer $n$). We will also assume that $S\bar\Theta$  has an increasing, unbounded argument (otherwise we replace $S$ with a large power). 
By Proposition 4.1 we have \begin{equation}\label{H2} N^+[U\Theta^{1-\epsilon}\bar\Theta]\ne0.\end{equation}Since the argument of $U\Theta^{1-\epsilon}$ is increasing, there is an  inner function $J$ such that $$U\Theta^{1-\epsilon}=XJ, \qquad \|\arg~X\|_\infty\le\pi.$$
From \eqref{H2} we have $N^+[J\bar\Theta]\ne0$, and so by Corollary 5.4
 \begin{equation}\label{H3}N^p[J\bar\Theta\bar S^{\epsilon}]\ne0,\qquad (p<1).\end{equation} Note that$$U\bar S^{2\epsilon}=\left(U\Theta^{1-\epsilon}\bar\Theta\bar S^{\epsilon}\right)~\left(\Theta^{\epsilon}\bar S^{\epsilon}\right)=X~\left(J\bar\Theta\bar S^{\epsilon}\right)~\left(\Theta^{\epsilon}\bar S^{\epsilon}\right).$$
Since the argument of $S^{\epsilon}\bar \Theta^{\epsilon}$ is increasing and unbounded,  we can find an infinite Blaschke product  $\Psi$ such that$$\Theta^{\epsilon}\bar S^{\epsilon}=Y\bar\Psi, \qquad  \|\arg~Y\|_\infty\le\pi.$$
Thus  the symbol $U\bar S^{2\epsilon}$ has the following representation:
$$U\bar S^{2\epsilon}=\left(J\bar\Theta\bar S^{\epsilon}\right)~(XY\bar\Psi), \qquad  \|\arg~XY\|_\infty\le2\pi,$$
and by \eqref{bmo2} we have
 \begin{equation}\label{H4}N^p[XY\bar\Psi]\ne0,\qquad (p<1/2).\end{equation}
Combining \eqref{H3} and \eqref{H4}, we get \eqref{H1} by H\"older's inequality.

\bs\no {\bf Remark.} It is clear from the proof that if $\gamma$ is almost decreasing, then   the \\$N^p$-kernels  are  infinite dimensional.

\bs\subsection{Proof of Theorem B} Recall that $J$ is a meromorphic  inner function, $S=e^{i\sigma}$ with $\sigma'(x)\asymp|x|^\kappa$, and $c=c(J,S;\kappa)$.
Applying Theorem A (or rather its corollary) to $U=J$ we conclude that if $a<c$ then $N^+[J\bar S^a]=0$ and therefore  $N^p[J\bar S^a]=0$ for all $p>0$. On the other hand, if $a>c$, then $N^p[J\bar S^a]\ne0$ for  some $p>0$, and in fact the kernel is infinite dimensional, as we just mentioned.  The following proposition completes the proof.  

\ms\no A unimodular function  $S$ is called {\it tempered} if $~\exists n$,\; $S'(x)=O\left(|x|^n\right)$ as $x\to\infty.$

\ms
\begin{prop}
 If $S$ is a tempered unimodular function, then for any meromorphic inner function $J$ and any $p>0$, 
$$\dim~N^p[J\bar S]=\infty\qquad\Rightarrow\qquad \dim~N^\infty[J\bar S]=\infty.$$\end{prop}

\ms\no \begin{proof} First of all we observe that the statement is true if $S$ is a tempered inner function, $S=\Theta$. By Carleson's type embedding theorem \cite{TV}, all elements in $N^p[J\bar \Theta]$ have at most polynomial growth at infinity, see details in \cite{MIF}, Section 4.1. Since the kernel is infinite dimensional, it contains functions with many zeros in $\C_+$. Dividing such functions by appropriate polynomials we obtain  functions in  $N^\infty[J\bar \Theta]$.

\ms\no Let now $S$ be an arbitrary tempered unimodular function. By Lemma 6.1 we can find a tempered inner function $\Theta$ and a bounded  real-valued function $\chi$ such that
$$S=\Theta\bar X,\qquad X=e^{2i\chi}.$$
By the previous observation, we have 
\begin{equation}\label{B1}\dim~N^\infty[J\bar \Theta]=\infty, \end{equation}
 and it remains to show that\begin{equation}\label{B2}\exists n,\quad N^\infty[X\bar b^n]\ne0,\end{equation}
($b$ is the Blaschke factor \eqref{bla}). Indeed, combining \eqref{B1} and \eqref{B2} we conclude that the kernel 
$$N^\infty[J\bar S\bar b^n]=N^\infty[J\bar \Theta ~X\bar b^n]$$ is infinite dimensional, which allows us to get rid of $b^n$.

\bs\no 
To prove \eqref{B2}, consider the outer function$$H=e^{\ti\chi-i\chi},\qquad{\rm so}\quad X=\frac{\bar H}H.$$
We will have $$(z+i)^{-n}H(z)\in N^\infty[X\bar b^n],\qquad (n\gg1)$$
if we can show that 
$h:=|H|=e^{\ti\chi}$ has at most polynomial growth at infinity.
Without loss of generality, we can assume that the $L^\infty$-norm of $\chi$ is so small that $h\in L^2_\Pi$. We have
\be\label{fi}h(x)-h(0)\le\int_0^x|h'|=\int_0^xh~|\ti\chi'|
\lesssim \|h\|_{L^2_\Pi}~(1+x^2)^{1/2}~\left(\int_0^x|\ti\chi'|^2\right)^{1/2}. \end{equation}
Since 
$|\chi'(t)|\lesssim|t|^n$  by construction,
for each $x>0$ we can represent $\chi$ as the sum of two smooth functions, $$\chi=\chi_1+\chi_2,$$ such that
$$\|\chi_1'\|_{L^2}\lesssim|x|^n,  \qquad \|\chi_2\|_{L^\infty} \asymp1,\qquad \chi_2=0 \quad{\rm on}\quad(-2x, 2x).$$
(For example, take $\chi_1=\phi\chi$, where $\phi$ is a smooth "bump" function such that $\phi$ is equal to $1$ on $(-2x, 2x)$ and $0$ on $\R\sm(-3x,3x)$.)
Then we have
$$\|\ti\chi_1'\|_{L^2}\lesssim|x|^n,  \qquad |\ti\chi'_2|\lesssim1 \quad{\rm on}\quad(0, x),$$
and so \eqref{fi} shows that $h$ has at most  polynomial  growth.  
\end{proof}

\bs

\end{document}